\input amstex

\def\b1{\text{\bf 1}}

\def\BC{{\Bbb C}}

\def\BQ{{\Bbb Q}}
\def\BZ{{\Bbb Z}}
\def\BR{{\Bbb R}}

\def\CL{{\Cal L}}

\def\cp{{\frak p}}


\def\hra{\hookrightarrow}
\def\iso{\buildrel\sim\over\longrightarrow}

\parskip=6pt

\documentstyle{amsppt}
\document
\magnification=1200
\NoBlackBoxes


\centerline{\bf  The Chiral de Rham Complex and}
\centerline{\bf Positivity
of the  Equivariant Signature of the Loop Space}

\bigskip

\centerline{ V. Gorbounov, F. Malikov\footnote{partially supported by an
NSF grant}}

\bigskip

{\bf Abstract}

 In this note we show that the positivity property of
the equivariant signature of the loop space, first observed in [MS1]
in the case of the even-dimensional projective spaces, is valid for
  Picard number 2 toric
varieties. A new formula for the equivariant
 signature of the loop space in the case of a toric spin variety is derived.

\bigskip\bigskip

{\bf 0. Introduction}
The equivariant signature of the loop space,
an example of an elliptic type genus, is, in particular, a rule which
assigns to a manifold $X$ a  power series in $q$:
$$
X\mapsto sign(q,\CL X)=b_{0}+b_{1}q+b_{1}q^{2}+\cdots .
$$
The intuition behind its definition (a well-known theorem in fact, see sect. 2) is
that $sign(q,\CL X)$ is the equivariant index of a certain elliptic operator
and therefore the coefficients $b_{j}$ are   integers,
 not necessarily
positive since they are equal to
 the difference between the dimensions of two vector
spaces. The cohomological computations of [MS1] fairly unexpectedly
 showed
 that in the case of the even-dimensional complex projective
space ${\Bbb P}^{2n}$ the series $sign(q,{\Cal L} {\Bbb P}^{2n})$ has positive
coefficients. Moreover, there is a graded vector space naturally
associated to ${\Bbb P}^{2n}$ so that $b_{j}$
essentially equals the dimension of the
$j$-th homogeneous component. Let us formulate this result more precisely.

Associated with any smooth manifold $X$ there is a sheaf of vertex algebras
$\Omega^{ch}_{X}$, the chiral de Rham complex [MSV]. It is graded:
$$
\Omega^{ch}_{X}=\oplus_{j\geq 0}(\Omega^{ch}_{X})_{j},
$$
and so are its cohomology groups. The following equality
of formal power series was proven in [MS1]
(cf. Theorem 4.1 below)
$$
sign(q,{\CL} {\Bbb P}^{2n})=
2\sum_{j=0}^{\infty} \text{dim}
 H^{0}(  {\Bbb P}^{2n},(\Omega^{ch}_{  {\Bbb P}^{2n}})_{j})q^{j} - 1.
\eqno{(0.1)}
$$
One can say that the vector space
 $H^{0}(  {\Bbb P}^{2n},\Omega^{ch}_{  {\Bbb P}^{2n}}) $, a vertex algebra in fact, provides a realization of
$sign(q,{\Cal L} {\Bbb P}^{2n})$.

The present note came out of the discussion of an earlier version of
[MS1], and (0.1) appearing in the final version of [MS1] is a result of this
discussion.  We further use the Borisov-Libgober formula for the equivariant
signature of the loop space in the case of a toric variety [BL] to extend the
positivity result to the Picard number 2 toric varieties.
(Recall that the projective spaces
exhaust the class of Picard number 1 smooth toric varieties.) The very nature
of this calculation does not allow us to conclude whether  there is
a graded vector space such that it is naturally associated to the manifold in question and
realizes this signature. The existence of such a vector space, be it
a vertex algebra or not, remains an open question.

A similar computation combined with Witten's rigidity theorem gives
 the following  formula for the equivariant signature
of the loop space in the case of a spin toric variety $X$
of complex dimension $d$ (cf. Theorem 6.2):
$$
sign(q,{\Cal L}X)=signX\frac{\epsilon^{-\frac{d}{4}}}{2^{d}},
\eqno{(0.2)}
$$
where $sign X$ is the signature of $X$ and $\epsilon$  a well-known modular form defined by (2.2) below.

The r.h.s. of this formula is not unfamiliar to vertex algebra
specialists. Let $V$ be  the vertex algebra of
$2d$
free bosons coupled to $2d$ free fermions. $V$ is naturally
graded:
$$
V=\oplus_{j=0}^{\infty} V_{j}
$$
so that
$$
\sum_{j=0}^{\infty} \text{dim}V_{j} q^{j}= \epsilon^{-\frac{d}{4}}.
\eqno{(0.3)}
$$
Having compared (0.1) to (0.2-0.3) one
  may perhaps be so bold as
to conjecture that the vertex algebra $V$ is isomorphic  to
$ H^{0}(  X,\Omega^{ch}_{ X })$. If correct this conjecture will be
 a natural
extension of the realization result (0.1) to the spin-case.

 {\bf Acknowledgments.} We are grateful to F.Hirzebruch and D.Zagier for
their interest in this work. Part of this work was completed when we were
visiting the Max-Planck-Institut f\"ur Mathematik in Bonn. We are grateful
to the Institut for the outstanding working conditions.

\bigskip\bigskip

{\bf 1. Genera.} In our brief review of the relevant genera we shall
follow the book [HBJ].

Let $R$ be a commutative ring, $\Omega^*$  the cobordism ring.
Recall that $\Omega^*\otimes \BQ=\BQ[P^2,P^4,...]$, where $P^{2n}$ is the cobordism class of the complex projective space of dimension $2n$. According to
 Hirzebruch, a {\it genus} is a ring homomorphism
 $$
 g:\;\Omega^*\otimes \BQ\rightarrow R.
 \eqno{(1.1)}
 $$
An invertible even formal power series $Q(x)=a_0 + a_1x^2 + \ldots$  with coefficients in $R$ defines a genus, to be denoted $g_{Q}$, as follows.
 Let $x_1,\ldots x_n$ be formal variables and $p_i$  the $i$-th elementary symmetric function in $x_i^2$. Then
$$Q(x_1)Q(x_2)\ldots Q(x_n)=a_0^n(1+K_1(p_1) + K_2(p_1,p_2) + \ldots)$$
where $K_i(p_1,\ldots ,p_i)$ is a (uniquely determined) homogeneous polynomial of degree $2i$.  For a manifold $X$ of dimension $4n$
define
$$
g_{Q}(X)=a_0^nK_n(p_1,\ldots p_n)[X],
\eqno{(1.2)}
$$
where $p_i(M)$, $i=1,...,n$, is the i-th the Pontryagin class of $X$ and $[X]$ is the fundamental class of $X$. In fact the rule
$Q(x)\rightarrow g_{Q}$ sets up a 1-1 correspondence.

\bigskip

{\bf 2. The elliptic genus and the equivariant signature of the loop space.}
 The notion of an elliptic genus is due to Ochanine[O].
A genus is called elliptic if it is associated to a series $Q(x)$
 such that:

if  $f(x)=\frac{x}{Q(x)}$, then
$$
(f^\prime)^2 = 1-2\delta f^2 + \epsilon f^4,
\eqno{(2.1)}
$$
for some parameters $\delta$ and $\epsilon$.

All solutions to (2.1) can be
 constructed
as follows. Fix a lattice $L$ in $\BC$ with generators
$\omega_{1},\omega_{2}$ and let $\cp_{L}(z)$ be
 the corresponding Weierstra\ss  function. Then
$$
f(x)=\frac{1}{ \sqrt{\cp_{L}(x)-\cp_{L}(\omega_{1}/2)} }
$$
is a solution to (2.1) with
$$
\delta=-\frac{2}{3} \cp_{L}(\omega_{1}/2) ,
$$
$$
\epsilon=[\cp_{L}(\omega_{1}/2) -\cp_{L}(\omega_{2}/2)]
[\cp_{L}(\omega_{1}/2) -\cp_{L}((\omega_{1}+\omega_{2})/2].
$$

Therefore each elliptic genus equals (cf. (1.2))
$$
 g_{x\sqrt{\cp_{L}(x)-\cp_{L}(\omega_{1}/2)}}(.),
$$
for some $L$ and $\omega_{1}$,
and $g_{x\sqrt{\cp_{L}(x)-\cp_{L}(\omega_{1}/2)}}(.)$, as a function
of $L$ and $\omega_{1}$, can be considered
 the universal elliptic genus.

Introducing the modular parameter $q=\exp{(2\pi i \omega_{2}/\omega_{1})}$
one sees that all the defined expressions are naturally identified
with functions of $q$. In particular, $\delta$ and $\epsilon$ become
the standard generators of the ring
$M_{*}(2)$ of modular
forms on $\Gamma_{0}(2)$ so that $\delta\in M_{2}(2)$,
$\epsilon\in M_{4}(2)$. Needed for (2.6) below is the following formula for
$\epsilon$:
$$
\epsilon=(2\prod_{n=1}^\infty\frac{(1+q^n)^2}{(1-q^n)^2})^{-4},
\eqno{(2.2)}
$$
Likewise, for any $X$,
 $g_{x\sqrt{\cp_{L}(x)-\cp_{L}(\omega_{1}/2)}}(X)$ is a function of
$q$. We shall emphasize this by changing (and unburdening) the notation
as follows:
$$
och(q,X)=g_{x\sqrt{\cp_{L}(x)-\cp_{L}(\omega_{1}/2)}}(X).
$$
In fact, for a  manifold $X$ of real dimension $4k$,
$och(q,X)$ is a weight $2k$ polynomial in $\delta$ counted with weight
2 and $\epsilon$  counted with weight 4. Therefore $och(q,X)$ is a
weight $2k$ modular
form on $\Gamma_{0}(2)$.

We will be more interested in a closely related genus,
 $sign(q,{\Cal L}X)$, proposed by Witten and called
the formal equivariant signature of the loop space [HBJ].
 By definition
$$
sign(q,{\Cal L}X)=g_{Q},
\eqno{(2.3a)}
$$
where
$$
Q(x)=x\frac{1+e^{-x}}{1-e^{-x}}\Pi_{n=1}^\infty \frac{(1+q^ne^{-x})(1+q^ne^{x})}{(1-q^ne^{-x})(1-q^ne^{x})}.
\eqno{(2.3b)}
$$

This genus can be equivalently defined to be the signature of $X$ twisted by the bundle
$$W=\otimes_{n=1}^{\infty} S_{q^n}TX_\BC\otimes_{n=1}^{\infty} \Lambda_{q^n}TX_\BC,
\eqno{(2.4)}
$$
where we habitually set
$$
\Lambda_{t}E = \sum_{i=0}^\infty t^i\Lambda^i E ,
\eqno{(2.5a)}
$$
$$S_{t}E = \sum_{i=0}^\infty t^iS^i E.
\eqno{(2.5b)}
 $$
 Yet another possibility is to define $sign(q,{\Cal L}X)$ to be the index of the elliptic operator $d+d^*$ acting on the global sections of the bundle
 (2.4).

We conclude this section by noting that $sign(q,{\Cal L}X)$ is
 related to the  elliptic genus as follows:
 $$
 sign(q,{\Cal L}X)=och(q,X)\epsilon^{-\frac{k}{2}},
 \eqno{(2.6)}
 $$
 where $\epsilon$ is defined by (2.2)
and $dim X=4k$.

\bigskip

{\bf 3. The chiral de Rham complex and the equivariant signature of the loop space.}   Defined in [MSV]
 for any smooth complex manifold $X$ of complex dimension $n$ there is a sheaf $\Omega^{ch}_{X}(X)$ of vertex algebras over $X$.
 Morally,  $\Omega^{ch}_{X}(X)$ is a semi-infinite de Rham complex on a ``small'' loop space with coefficients
  in distributions supported on the submanifold of
 analytically contractible loops. This vague assertion was made a theorem in [KV] after overcoming considerable
 technical difficulties.

$\Omega^{ch}_{X}(X)$ is not a sheaf of ${\Cal O}_X$-modules but it possesses a filtration such that the associated graded sheaf $gr \Omega^{ch}_{X}$ is. In fact,
$gr \Omega^{ch}_{X}$ is associated to a holomorphic vector bundle and
 this vector bundle was explicitly described in [BL]
 as follows:
$$
gr \Omega^{ch}_{X} =\otimes_{n=1}^{\infty} \{S_{q^n}TX\otimes S_{q^n}TX^*\otimes \Lambda_{yq^{n-1}}T^*X\otimes\Lambda_{y^{-1}q^{n}}TX\},
\eqno{(3.1)}
$$
where we use the notation introduced in (2.5a,b). Thus $gr \Omega^{ch}_{X}(X)$ is bi-graded so that the component of weight $(i,j)$ is the coefficient
of $q^{j}y^{i}$. In fact, the sheaf $\Omega^{ch}_{X}$ is itself bi-graded :
$$
\Omega^{ch}_{X}=\oplus_{i\in\BZ,j\geq 0}(\Omega^{ch}_{X})^{i}_{j},
\eqno{(3.2)}
$$
and this bi-grading descends to the graded object. In particular, $\Omega_{X}^{i}$, i.e. the sheaf of
holomorphic $i$-forms , canonically identifies
with $(\Omega^{ch}_{X})^{i}_{0}$:
$$
 \Omega_{X}^{i}\iso (\Omega^{ch}_{X})^{i}_{0},
\eqno{(3.3.)}
$$
which partially justifies the name `` a chiral de Rham complex''.

 Consider the Euler character $\text{Eu} (\Omega^{ch}_{X})(q,y)$, which by definition is a formal Laurent power series in $q,y$ such that the coefficient
 of $y^{i}q^{j}$ is the Euler characteristic of the component $(\Omega^{ch}_{X})^{i}_{j}$.  Define $\text{Eu} (\Omega^{ch}_{X})(q)$ to be
 $\text{Eu} (\Omega^{ch}_{X})(q,1)$
Applications of $\Omega^{ch}_{X}$ to elliptic genera are based on the following observation due to [BL]:
$$
\text{Eu} (\Omega^{ch}_{X})(q)=sign(q,{\Cal L}X).
\eqno{(3.4)}
$$

 Let us prove (3.4) for the sake of completeness. Since the Euler characteristic does not change under the passage to the graded object, we can write

$$
\text{Eu} (\Omega^{ch}_{X})(q)=\text{Eu} (gr \Omega^{ch}_{X})(q)=
\int_X ch(gr\Omega^{ch}_{X}) tdX,
\eqno{(3.5)}
$$
where the 2nd equality follows from the Riemann-Roch Theorem, and $tdX$ is the Todd genus of $X$. We now compute $ch(gr\Omega^{ch}_{X})$ by using (3.1)
and the multiplicativity of $ch$ to the effect that
$$
ch(gr \Omega^{ch}_{X})=\prod_{n=1}^{\infty} \{ch(S_{q^n}TX)
 ch (S_{q^n}TX^*)ch(\Lambda_{q^{n-1}}T^*X) ch(\Lambda_{q^{n}}TX)\}.
\eqno{(3.6)}
$$

 As is well known for a $k$ dimensional manifold $X$
$$
ch \Lambda_{q^{n-1}}TX = \Pi_{i=1}^{k} (1+q^{n-1}e^{x_i}),\;\;\;ch \Lambda_{q^{n-1}}TX^{*} = \Pi_{i=1}^{k} (1+q^{n-1}e^{-x_i}),
$$
$$
ch S_{q^{n}}TX = \Pi_{i=1}^k \frac{1}{(1-q^{n}e^{x_i})},\;\;\;ch S_{q^{n}}TX^{*} = \Pi_{i=1}^k \frac{1}{(1-q^{n}e^{-x_i})},
$$
$$
tdX=\Pi_{i=1}^{k}\frac{x_{i}}{1-e^{-x_{i}}}.
$$
Plugging  these in the right hand side of
(3.6) we readily see
   that the integrand in the right hand side of (3.5) is
   $Q(x_{1})Q(x_{2})\cdots Q(x_{k})$, where $Q(x)$ is the series (2.3b).
    Therefore (3.4) is  identical to the definition (1.2) with
 $Q(x)$ given by (2.3b).
$\qed$

Note that $\text{Eu} (\Omega^{ch}_{X})(0)$ is equal to the signature of the manifold $X$.

\bigskip\bigskip

{\bf 4. Positivity of the equivariant signature of the loop space.}

Introduce
the formal character
$$
ch H^{i}(X,\Omega^{ch}_{X})=\sum_{j=0}^{\infty}\text{dim}H^{i}(X,(\Omega^{ch}_{X})_{j})q^{j},
$$
where we  ignore the upper-index grading, cf. (3.2). (It follows easily from (3.1) that the dimensions entering the formal character are all finite.)

{\bf Theorem 4.1.} The following character formula is valid
$$
2 ch H^{0}(  {\Bbb P}^{2n},\Omega^{ch}_{  {\Bbb P}^{2n}}) - 1=sign(q,{\Cal L} {\Bbb P}^{2n}).
$$
\bigskip

{\bf Corollary 4.2.}
The series $sign(q,{\Cal L} {\Bbb P}^{2n})$ has positive coefficients and the coefficients against positive powers of $q$ are even.

\bigskip

Corollary 4.2 is an immediate consequence of Theorem 4.1 whereas Theorem 4.1 is a simple consequence of certain properties of the chiral de Rham complex
over projective spaces discovered in [MS1]. We shall list these properties here and reproduce from [MS1]  a simple computation  leading to Theorem 4.1.

 The definition of the Euler characteristic and (3.5) give
 $$
 sign(q,{\Cal L} {\Bbb P}^{2n})=\sum_{i=0}^{2n}(-1)^{i}
 ch H^{i}(  {\Bbb P}^{2n},\Omega^{ch}_{  {\Bbb P}^{2n}}).
 \eqno{(4.1)}
 $$

We now make use of the following

{\bf Theorem [MS1].} The natural embedding of sheaves $\Omega^*_{{\Bbb P}^n} \hra \Omega^{ch}_{  {\Bbb P}^{n}}$ due to (3.3) provides an isomorphism
$$H^i({\Bbb P}^n,\Omega^*_{{\Bbb P}^n})\iso  H^i(  {\Bbb P}^{n},\Omega^{ch}_{  {\Bbb P}^{n}}). $$
for $0<i<n$.

\bigskip

This assertion reduces (4.1) to
$$
 sign(q,{\Cal L} {\Bbb P}^{2n})=ch H^{0}(  {\Bbb P}^{2n},\Omega^{ch}_{  {\Bbb P}^{2n}})+
 ch H^{2n}(  {\Bbb P}^{2n},\Omega^{ch}_{  {\Bbb P}^{2n}})-1.
 \eqno{(4.2)}
 $$

 To conclude it remains to notice that due to the chiral Poincar\'e duality [MS2]
 $$
 ch H^{0}(  {\Bbb P}^{2n},\Omega^{ch}_{  {\Bbb P}^{2n}})=
 ch H^{2n}(  {\Bbb P}^{2n},\Omega^{ch}_{  {\Bbb P}^{2n}}).
 $$
 $
 \qed
 $

{\bf Remark} The vector space $H^0({\Bbb P}^n,\Omega^{ch}_{{\Bbb P}^n})$ is a vertex algebra. Thanks to Theorem 4.1, the known [HBJ] modular
properties of the equivariant signature of the loop space say that $ch H^0({\Bbb P}^n,\Omega^{ch}_{{\Bbb P}^n})$
  is a modular function when $n=0$ mod $4$. The modular properties of  characters  have been the hallmark of vertex algebra theory,
 see for example [Z].

{\bf 5. Extending the positivity result}

{\bf 5.1} It seems natural to ask in what generality the positivity
result (Corollary 4.2) holds true. Since according to [BR] every cobordism class contains a non-singular toric variety, not
every non-singular toric variety $X$ has positive $sign(q,\CL X)$. We will show nevertheless that apart from projective spaces there is a class of toric varieties with this
property. To do so, we will  have to rely  on the calculations from [BL]
instead of sheaves of vertex algebras.

When talking about toric varieties we shall keep to the
 following notation. Let
  $e_i$, $1\leq i\leq d$, be
 the standard basis of $\BZ^d$; thus the j-th component of $e_i$ is
$\delta_{ij}$. Define the inner product
$$
\BZ^{d}\times\BZ^{d}\rightarrow \BZ,\; x,y\mapsto x\cdot y
$$
by the requirement that
$$
e_{i}\cdot e_{j}=\delta_{ij}.
$$
This identifies $\BZ^{d}$ and hence $\BR^{d}$ with their duals.

By $\Sigma$ we denote a {\it complete, regular} fan in $\BR^{d}$.
This means, in particular, that each cone $C^{*}\in\Sigma$ is spanned
by part of a basis of $\BZ^{d}$, and we denote this spanning
set by $|C^{*}|$.

Associated to  $\Sigma$ there is a smooth compact toric variety
of complex dimension $d$ to be denoted $X_{\Sigma}$ or simply $X$
if no confusion is likely to arise.

  The following formula holds
true [BL]
$$
sign(q,{\Cal L}X)=\sum_{m\in \BZ^{d}}\sum_{C^*\in \Sigma}
(-1)^{codim C^*}(\prod_{n\in |C^*|}\frac{1}{1+q^{m\cdot n}})\epsilon^{-d/4}.
\eqno{(5.1)}
$$
 Note that
to make sense out of this expression one has to expand each factor
$$
\frac{1}{1+q^{m\cdot n}}
$$
at $q=0$ and then  convince oneself that the sum of thus arising power
series with respect to $m\in \BZ^{d}$ makes sense as a formal power series.

{\bf Remark} The formula in  Theorem 5.5 in [BL]   where we
borrowed (5.1) contains an extra factor $(-1)^{d/2}$. We drop
it so  as to conform to the standard notation.
\bigskip

{\bf 5.2}
Now we extend the result of  Theorem 4.1 to a larger class of
toric variates.
Recall that the Picard number of a toric variety is the difference between the dimension of the variety and the number of one dimensional cones. Each smooth
Picard number one toric variety is a  projective space. The
Picard number two
toric varieties, which can be viewed
as generalized Hirzebruch's surfaces, were classified in [Kl].  Let us formulate
this result.

We give ourselves a triple of integers
 $d$, $s$, $r$  such that $1 < d$, $1 < s < d+1$, $r = d-s+1$,
and an increasing sequence of non-negative integers  $a_1,\ldots a_r$.   Define the following vectors in $\BZ^{d}$:

$$u_i = e_i, 0<i<r+1;$$
$$u_{r+1} = -\sum_{i=1}^r u_i;$$
$$v_j = e_{r+j}, 0<j<s;$$
$$v_s = \sum_{i=1}^r a_ie_i -\sum_{j=1}^{s-1}v_j;$$

We set $U=\{u_1,\ldots , u_{r+1}\}$, $V = \{v_1,\ldots ,v_s\}$ and
let $C^{*}_{ij}\subset \BR^{d}$, $0< i < r+2$, $0< j<s+1$,  be the cone
spanned
by  $U\cup V \setminus \{u_i,v_j \})$. One checks that there is a
uniquely determined regular complete fan such that
$\{C^{*}_{ij},\;0< i < r+2,0< j<s+1\}$
is the set of $d$-dimensional cones. We denote this fan by
$\sum_d(a_1,\ldots a_r)$ and the corresponding toric variety by
$X_d(a_1,\ldots a_r)$.

{\bf Theorem 5.1} [Kl] Every compact smooth toric
variety of complex dimension $d$ with $d+2$ generators is isomorphic to precisely one of the varieties $X_d(a_1,\ldots a_r)$.

{\bf Theorem 5.2}

1. For $d$ even

$$
sign(q,{\Cal L} {\Bbb P^d})=\displaystyle{\sum_{m\in \BZ^{d}}\frac{2}{(1+q^{-m_1\ldots -m_d})\prod^{d}_{i=1}(1+q^{m_i})}\epsilon^{-d/4}},
\eqno{(5.2)}
$$

where $m=(m_1\ldots m_d)$. The series in the RHS has positive
coefficients.

2. For a smooth toric variety $X$ of even complex dimension $d$
with  Picard number $2$, $sign(q,{\Cal L} X)$ is 0 if $s$ is even.
Otherwise
$$
sign(q,{\Cal L} X)=\sum_{m\in \BZ^{d}}\frac{2+2q^{m\cdot (v_1+\ldots v_s)}}{\prod_{i=1}^{r+1} (1+q^{m\cdot u_i}) \prod_{j=1}^s (1+q^{m\cdot v_j})}\epsilon^{-d/4}
\eqno{(5.3)}
$$
The series in the RHS has positive
coefficients.

\bigskip

{\bf Remark.} Statement 1. of Theorem 5.2 is an alternative to the
chiral de Rham complex approach of sect. 4. But note that
formula (5.2) is  of a different nature than that in Theorem 4.1,
and the comparison of (5.2) and Theorem 4.1 may give rise to non-trivial
combinatorial identities.

\bigskip

{\bf Proof.} Both of the statements follow from (5.1) and
 the explicit description of the
fans of the toric varieties in question.

1.
Let $k_1=e_1, \ldots ,k_d=e_d, k_{d+1}=-e_1\ldots -e_d$. The fan
defining $\Bbb P^d$
 consists of the cones spanned by all proper subsets of the set
 $\{k_1, \ldots ,k_{d+1}\}$.

Now observe that (5.1) rewrites as follows:
$$
sign(q,{\Cal L}X)=\epsilon^{-d/4}\sum_{m\in M}<m, \sum_{C^*\in \Sigma}(-1)^{codim C^*}(\prod_{k\in|C^*|}\frac{1}{1+q^{k}})>,
\eqno{(5.4)}
$$
where
$q^{k}$ is an element of the group ring of $\BZ^{d}$  and
$<m,q^{k}>=q^{m\cdot k}$.

The expression
$$
 \sum_{C^*\in \Sigma}(-1)^{codim C^*}(\prod_{k\in|C^*|}
\frac{1}{1+q^{k}}),
\eqno{(5.5)}
$$
(appearing in the  r.h.s. of (5.4)) in the case of the projective space
simplifies as follows:
$$
\sum_{B}
(-1)^{d-\# B}\frac{1}{\prod_{k\in B}(1+q^{k})},
\eqno{(5.6)}
$$
where the summation is performed over all   {\it proper} subsets
 $B\subsetneq\{k_1,\ldots ,k_{d+1}\}$ and $\# B$ is the number of elements in $B$.
Converting (5.6)  to the common denominator we obtain
$$
\frac{\sum_{B}(-1)^{d-\# B}\prod_{k\in \bar B}(1+q^{k}) }
{\prod^{d+1}_{i=1}(1+q^{k_{i}})},
\eqno{(5.7)}
$$
where $\bar B$ denotes the complement of  $B$.
The multiple use of the binomial identuty
$$
\sum_{i=0}^{n}(-1)^{i}C_{n}^{i}=0\text{ if } n>0
\eqno{(5.8)}
$$
allows us to collect the like terms in the numerator of (5.7). The result
is
$$
\frac{1+(-1)^{d}q^{k_1+\ldots +k_{n+1}}}{\prod^{d+1}_{i=1}(1+q^{k_i})}.
\eqno{(5.9)}
$$
Since $k_{1}+k_{2}+\cdots +k_{d+1}=0$, it
 is zero if $d$ is odd  (as it should !) and
$$
\frac{2}{\prod^{d+1}_{i=1}(1+q^{k_i})}
\eqno{(5.10)}
$$
otherwise. Plugging (5.10) in the r.h.s. of (5.4) and
performing  pairing with $m\in\BZ^{d}$ we obtain (5.2),
as desired.
It is clear that every term in the denominator of (5.2) cancels against the appropriate term in $\epsilon^{-d/4}$.
Indeed, for a fixed $m\in\BZ^{d}$ the denominator of (5.2) contains at most $d+1$
factors  $(1+q^n)$ for each $n\in \Bbb Z$ whereas
$\epsilon^{-d/4}$, see (2.2), contains $2d$ such factors
 in the numerator. Therefore all the coefficients in the
series $sign(q,{\Cal L} {\Bbb P^d})$ are positive.

2. Proof of the second statement is similar.
 Theorem 5.1 combined with (5.1) gives the following analogue of (5.4):

 $$
sign(q,{\Cal L} X)=\epsilon^{-d/4}\sum_{m\in M}<m,\sum_{I,J}\frac{(-1)^{d-\# I - \# J}}{\prod_{i\in I}(1+q^{u_i})\prod_{j\in J}(1+q^{v_j})}>,
\eqno{(5.11)}
$$
where $I$ and $J$ are {\it proper} subsets of $\{1,\ldots ,r+1\}$ and $\{1,\ldots ,s\}$ respectively and the vectors $u_i$, $v_j$ are
those defined in the beginning of 5.2.

Converting  the sum $\sum_{I,J}$ in (5.11)
to the common
denominator we obtain

$$
\sum_{I,J}\frac{(-1)^{d-\# I - \# J}\prod_{i\in \bar I}(1+q^{u_i})\prod_{j\in \bar J}(1+q^{v_j})}
{\prod_{i=1}^{r+1}(1+q^{u_i})\prod_{j=1}^{s}(1+q^{v_j})}
\eqno{(5.12)}
$$

Now observe that the numerator of (5.12) is the product of two factors
analogous to the numerator of (5.7) -- one is  the numerator of (5.7)
with $B$ replaced with $\{u_{i},\, i\in I\}$, another is also
with $B$ replaced with $\{v_{j},\, j\in J\}$.
Therefore identity (5.8), which allowed us to pass from (5.7) to (5.9),
allows us to analogously rewrite (5.12) as follows:
$$
\frac{(1+(-1)^{r}q^{u_{1}+u_{2}+\cdots +u_{r+1}})(1+(-1)^{s-1}q^{v_{1}+v_{2}
+\cdots +v_{s}})}
{\prod_{i=1}^{r+1}(1+q^{u_i})\prod_{j=1}^{s}(1+q^{v_j})}.
\eqno{(5.13)}
$$
Since by definition (see the beginning of 5.1)
$$
\sum_{i=1}^{r+1}u_i=0,
$$
 (5.13) vanishes if $r$ is odd and equals
$$
\frac{2(1+q^{v_{1}+v_{2}
+\cdots +v_{s}})}
{\prod_{i=1}^{r+1}(1+q^{u_i})\prod_{j=1}^{s}(1+q^{v_j})}
\eqno{(5.14)}
$$
otherwise. (Note that $(-1)^{s-1}$ has disappeared because, $d$ being even, the relation
$d=r+s-1$ forces $r$ and $s$ to have different parity.)

Plugging (5.14) in the r.h.s. of (5.11) and performing  pairing with
$m\in\BZ^{d}$ we obtain desired formula (5.3).

To show that  series (5.3)
 has positive coefficients we  observe that for a fixed $m\in\BZ^{d}$
the denominator of (5.3) contains at most $r+s+1$, that is, $d+2$
factors  $(1+q^n)$ for each $n\in \Bbb Z$ whereas
$\epsilon^{-d/4}$, see (2.2), contains $2d$ such factors
 in the numerator. Therefore, having carried out the cancellations we make
(5.3) into a sum of power series with positive coefficients.
 $\qed$

\bigskip

{\bf 6. Toric spin varieties}

Recall Witten's rigidity theorem [W] proved in [BT].
 Let a torus $T^n$ act on a manifold X. This action
lifts to an action on the holomorphic bundle $gr \Omega^{ch}_{X}$,
see (3.1) Therefore each cohomology group $H^{i}(X,gr \Omega^{ch}_{X})$
becomes a $T^n$-module -- a direct sum of the torus characters
$t\mapsto t^m$, $t\in T^n$,
$m\in \Bbb Z^n$, in fact. Formula (3.4) then implies that
$sign(q,{\Cal L}X)$ is a formal sum of the torus characters and one can think
of $sign(q,{\Cal L}X)$ as a function of $t\in T^n$ with values in
$\BC[[q]]$.

{\bf Theorem 6.1 [BT]} If $X$ is a spin manifold equipped with an
 action of a torus $T^n$, then
$sign(q,{\Cal L}X)$ is a constant function of $t\in T^n$.

{\bf Theorem 6.2} For a toric spin variety $X$ of complex dimension
 $d$,
$$
sign(q,{\Cal L}X)=signX\frac{\epsilon^{-\frac{d}{4}}}{2^{d}}.
$$
\bigskip

{\bf Proof.} Any toric variety $X$ carries the natural action of
 a torus;
hence $sign(q,{\Cal L}X)$ is a formal sum of
the torus characters. Formula (5.1) sharpens accordingly [BL]:
$$
sign(q,{\Cal L}X)=\sum_{m\in M}t^m\sum_{C^*\in \Sigma}(-1)^{codim C^*}(\prod_{i=1,\ldots dim C^*}\frac{1}{1+q^{m\cdot n_i}})\epsilon^{-d/4}.
\eqno{(6.1)}
$$
Theorem 6.1 implies that  only $t^0$ may appear in the r.h.s. of
(6.1) with non-zero coefficient. Therefore
$$
sign(q,{\Cal L}X) = C\epsilon^{-d/4}
\eqno{(6.2)}
$$
for some constatnt $C$. To compute $C$ recall that
$$
sign(0,{\Cal L}X)=sign X.
$$
Therefore, having specialized (6.2) to $q=0$ we obtain
$$
sign X= C 2^{d},
$$
and Theorem 6.2 follows. $\qed$

\bigskip
{\bf Remarks.}

1) According to [HS]  $sign(q,{\Cal L}X)=signX$ if $X$ is a
homogeneous space and a spin manifold at the same time.
Therefore Theorem 6.2  is
an extension of this result to spin toric varieties.

2) Theorem 6.2 suggests that perhaps the vertex algebra of $2d$ bosons
coupled to $2d$ fermions provides a natural realization of
$sign(q,{\Cal L}X)$  in the case of a toric spin manifold as discussed
in greater detail in the introduction .




\centerline{\bf References}

\bigskip\bigskip


[BL] L.~Borisov, A.~Libgober, Elliptic Genera and Applications to Mirror Symmetry,  Inv. Math. {\bf 140} (2000), 453-485;.

[BT] R. Bott, C. Taubes, On the rigidity theorems of Witten. J. Amer. Math. Soc. 2 (1989), no. 1, 137--186.

[D] A. N. Dessai, Elliptic genera, positive curvature and symmetry, preprint, 2002.

[GMS] V. Gorbounov, F. Malikov, V. Schechtman, Gerbes of chiral differential operators {\bf II} and {\bf III} math AG/0003170, math AG/0005201.

[HBJ] F. Hirzebruch, Th. Berger, R. Jung, Manifolds and modular forms. With appendices by Nils-Peter Skoruppa and by Paul Baum. Aspects of Mathematics, E20. Friedr. Vieweg $\&$
Sohn, Braunschweig, 1992. xii+211 pp.

[HS] F. Hirzebruch, P. Slodowy, Elliptic genera, involutions, and homogeneous spin manifolds. Geom. Dedicata 35 (1990), no. 1-3, 309--343.

[K] V.~Kac, Vertex algebras for beginners,
Second Edition, University Lecture Series, {\bf 10},
American Mathematical Society, Providence, Rhode Island, 1998.

[Kl] P. Kleinschidt, A classification of toric varieties with few generators, Aerquationes Mathematicae 35 (1988) 254-266.

[KV] M. Kapranov, E. Vasserot,  Vertex algebras and the formal loop space, math.AG/0107143.

[MSV] F.~Malikov, V.~Schechtman, A.~Vaintrob, Chiral de Rham complex,
{\it Comm. Math. Phys.}, {\bf 204} (1999), 439-473.

[MS1] F.~Malikov, V.~Schechtman, Deformations of chiral algebras and quantum cohomology of toric varieties, to appear in Comm. Math. Phys.

[MS2] F.~Malikov, V.~Schechtman, Chiral Poincar\'e duality,
   Math. Res. Lett.  vol. {\bf 6} (1999), 533-546.

[O] S.~Ochanine, Sur les genres multiplicatifs d\'efinis par des int\'egrales elliptiques. (French) [On multiplicative genera defined by elliptic integrals] Topology 26 (1987), no. 2, 143--151.

[W] E.~Witten, The index of the Dirac operator in loop space,
in: Elliptic curves and modular forms in algebraic topology
(Princeton, NJ, 1986), 161-181, Lect. Notes in Math. {\bf 1326} (1988).

[Z] Y. Zhu, Modular invariance of characters of vertex operators algebras, J. Amer. Math. Soc. {\bf 9} (1996), no 1, 237--302.

\bigskip

\enddocument